# TRANSFORM MARTINGALE ESTIMATING FUNCTIONS

By T. Merkouris

*Statistics Canada*

An estimation method is proposed for a wide variety of discrete time stochastic processes that have an intractable likelihood function but are otherwise conveniently specified by an integral transform such as the characteristic function, the Laplace transform or the probability generating function. This method involves the construction of classes of transform-based martingale estimating functions that fit into the general framework of quasi-likelihood. In the parametric setting of a discrete time stochastic process, we obtain transform quasi-score functions by projecting the unavailable score function onto the special linear spaces formed by these classes. The specification of the process by any of the main integral transforms makes possible an arbitrarily close approximation of the score function in an infinite-dimensional Hilbert space by optimally combining transform martingale quasi-score functions. It also allows an extension of the domain of application of quasi-likelihood methodology to processes with infinite conditional second moment.

**1. Introduction.** Maximum likelihood estimation of parameters of discrete time stochastic processes is often not feasible because an explicit expression for the associated likelihood function is either unavailable or too complicated. In a wide variety of situations, however, a description of the process by an integral transform, such as the conditional characteristic function or the conditional Laplace transform, is more readily available than explicit likelihood or score functions. A broad array of such processes encountered in the literature includes the following. First, there are linear processes with infinite variance used in modeling certain time series phenomena. In particular, models in economics and in signal processing involving linear time series with error having a stable distribution have been considered in the literature; see, for example, McCulloch [25] and Nikias and Shao [30].

In general, no closed form expression suitable for likelihood methods exists for the densities of such processes, but a simple form of characteristic function or Laplace transform is available. Second, there are two rich classes of discrete and non-Gaussian continuous variate time series models. One of these classes includes models with exponential or gamma marginal distributions, which are of importance in queuing and network processes. The other class includes stationary processes with discrete marginal distributions, such as Poisson, geometric or negative binomial, which are useful for modeling counting processes consisting of dependent random variables. A detailed survey of first-order autoregressive processes with such distributions is given in Grunwald, Hyndman, Tedesco and Tweedie [17]. For the higher-order autoregressive case, see Billard and Mohamed [3], Alzaid and Al-Osh [2] and references therein. These non-Gaussian models have an intractable likelihood function because of its complexity or inherent discontinuities. Nonetheless, they are described handily by their Laplace transforms or probability generating functions. Third, there are aggregate models, such as aggregate Markov chains and the related compartmental models used in many areas, such as economics, population theory and the social sciences, when only aggregate data are available; see, for example, McLeish [26] and Leitnaker [24]. Such models involve convolutions whose densities are rarely tractable but have a simple representation by their probability generating functions.

In all these examples maximum likelihood estimation for the parameters of the relevant models has been regarded as unworkable, and other methods of estimation, typically using only the first two moments of the underlying distribution, are generally suboptimal.

In this article we propose an estimation method for discrete time stochastic processes that are conveniently specified by conditional integral transforms. The development of this method follows earlier work (Merkouris [27, 28]) that introduced quasi-likelihood estimation for discrete time semimartingales based on conditional integral transforms. Though quite distinct, this approach is related to previous estimation methods in the literature that involved fitting an empirical transform to its theoretical counterpart in the setting of i.i.d. random variables. In this sense, it is akin to the generalized moment procedure proposed by Feuerverger and McDunnough [14] and Brant [4] as an approximate maximum likelihood procedure.

In the parametric setting of a discrete time stochastic process with intractable likelihood function we build classes of martingale estimating functions by means of an integral transform that specifies the process. Such classes of transform-based martingale estimating functions fit into the general quasi-likelihood framework given by Godambe and Heyde [15]. In our quasi-likelihood approach, we obtain quasi-score functions as projections of the unavailable score function onto the special linear spaces formed by these classes. Thus, these transform quasi-score functions provide best linear approximations to the score function in Hilbert space. In contrast to the semiparametric setting of the ordinary quasi-likelihood, the specification of the process by its transform structure enables the utilization of distributional information beyond the second-order moment structure. Enlarged classes of suitable composite transform martingales can then be constructed, which may lead to an arbitrarily close approximation of the score function. The transform structure allows also an extension of the quasi-likelihood methodology to processes with infinite conditional second moment. Furthermore, combining transform martingale estimating functions may be effective in dealing with problems of identifiability of vector parameters.

The proposed transform martingale estimating functions are generally nonlinear in the observations, and for the main integral transforms they may be represented as perturbed polynomial quasi-score functions. In particular, a basic transform martingale estimating function can be expressed as a perturbed ordinary quasi-score (weighted conditional least squares) estimating function.

Other estimation methods based on empirical transforms for discrete time stochastic processes have appeared in the literature. Feuerverger [11] discussed an asymptotically efficient estimation procedure based on the "polycharacteristic" function in the setting of univariate stationary time series models. Brockwell and Liu [5] used the empirical characteristic function in estimation for a linear process with stable innovations. Estimators based on the Laplace transform were used for the estimation of parameters of progression time distributions in multi-stage models by Schuh and Tweedie [32], Feigin, Tweedie and Belyea [10] and Hoeting, Tweedie and Olver [20]. Abraham and Balakrishna [1] used the empirical Laplace transform to estimate a parameter of a first-order inverse Gaussian autoregressive process. Yao and Morgan [37] used a least squares approach based on empirical transforms for a class of indexed stochastic models. Each of the aforementioned methods is essentially an ad hoc approach to the particular estimation problem. In contrast, this article presents a general estimation procedure that is statistically and computationally efficient, of broad applicability and based on a comprehensive estimation theory.

The article is organized as follows. In Section 2 the transform method is introduced as a special quasi-likelihood estimation procedure with a potential of high efficiency that rests on the capacity of constructing combinations of transform martingale estimating functions. A link is established between the method of Feuerverger and McDunnough, and Brant, and the quasi-likelihood approach through the concept of best linear approximation of the score function that underlies both estimation procedures. In Section 3 the formulation of orthogonal projection of the score function onto suitable infinite-dimensional spaces lays the ground for the construction of potentially fully efficient transform martingale estimating functions. A practicable procedure for forming composite transform martingale estimating functions of nondecreasing efficiency is then developed. Computational issues related to the proposed estimating procedure are discussed in Section 4. Comparisons of martingale estimating functions based on important transforms are presented in Section 5. Also in Section 5, the relationship of the transform estimation method with ordinary quasi-likelihood and conditional least squares is established, and suggestions as to the choice of an efficient transform estimating function are made. Two brief examples illustrating important features of the transform method are given in Section 6.

## 2. Transform martingale estimating functions.

2.1. *The basic structure.* Let $\{\mathbf{Y}_j, 1 \le j \le n\}$ be a sample from a discrete time stochastic process which takes values in $r$-dimensional Euclidean space and whose distribution depends on a parameter $\boldsymbol{\theta}$ belonging to an open subset $\boldsymbol{\Theta}$ of $p$-dimensional Euclidean space. Suppose that the possible probability measures for $\{\mathbf{Y}_j\}$ are $\{P_{\boldsymbol{\theta}}, \boldsymbol{\theta} \in \boldsymbol{\Theta}\}$ and that each $(\Omega, \mathcal{F}, P_{\boldsymbol{\theta}})$ is a complete probability space. Let $\mathcal{F}_j$ denote the past-history sub-$\sigma$-field of $\mathcal{F}$ generated by $\mathbf{Y}_1, \ldots, \mathbf{Y}_j, j \ge 1$. Suppose that the conditional density function $f_{\boldsymbol{\theta}}(\mathbf{Y}_j | \mathbf{Y}_1, \ldots, \mathbf{Y}_{j-1})$ and the gradient of $\log f_{\boldsymbol{\theta}}(\mathbf{Y}_j | \mathbf{Y}_1, \ldots, \mathbf{Y}_{j-1})$, denoted by $\mathbf{s}_j$, exist. Allowing differentiation under an integral sign, $E(\mathbf{s}_j | \mathcal{F}_{j-1}) = \mathbf{0}$ almost surely (a.s.) and, thus, the score function $\mathbf{S}_n = \sum_{j=1}^{n} \mathbf{s}_j$ defines a zero-mean martingale, $\{\mathbf{S}_n, \mathcal{F}_n\}$, which is assumed to be square integrable.

The score function may be unavailable or too difficult to compute, so that estimating $\boldsymbol{\theta}$ by maximum likelihood is not feasible. Nevertheless, a class of workable martingale estimating functions may instead be constructed based on an integral transform.

We assume at first, for simplicity of exposition, that we deal with real valued random variables $Y_1, \ldots, Y_n$ whose distribution depends on a scalar parameter $\theta \in \Theta$. Then for the $j$th time point we write

$$F_j(y|\mathcal{F}_{j-1}) = P(Y_j \le y|\mathcal{F}_{j-1}), \qquad \hat{F}_j(y) = I_{[Y_j \le y]}, \qquad 1 \le j \le n,$$

where $I$ denotes the indicator function. For an indexed set of real or complex valued functions $\{g_t(Y), t \in T \subseteq \mathbb{R}\}$, the kernel class, we consider the integral transform

$$(1) \qquad c_j(t) = \int g_t(y) \, dF_j(y|\mathcal{F}_{j-1}) = E(g_t(Y_j)|\mathcal{F}_{j-1}),$$

where the kernel $g_t(\cdot)$ is such that the integral exists and is finite for all $\theta \in \Theta$ and all $t \in T$. The dependence of $F_j(y|\mathcal{F}_{j-1})$ and $c_j(t)$ on $\theta$ is suppressed notationally for convenience of writing. The important transforms include the characteristic function, the moment generating function and the probability generating function, with associated kernel classes $\{\sin(tY), \cos(tY), t \in \mathbb{R}\}$, $\{\exp(tY), t \in \mathbb{R}\}$ and $\{t^Y, t \in \mathbb{R}\}$, respectively, though others (e.g., the Laplace transform or the sequence of moments) may be used as befits the context.

Now let

$$\hat{c}_j(t) = \int g_t(y) \, d\hat{F}_j(y) = g_t(Y_j) \tag{2}$$

and write

$$h_j(t) = \hat{c}_j(t) - c_j(t) = g_t(Y_j) - E(g_t(Y_j)|\mathcal{F}_{j-1}), \qquad 1 \le j \le n. \tag{3}$$

Since $E(h_j(t)|\mathcal{F}_{j-1}) = 0$ a.s., for fixed $t \in T$ the $\{h_j(t), \mathcal{F}_j\}$ are martingale differences of a zero-mean martingale, say, $\{H_n(t) = \sum_{j=1}^n h_j(t), \mathcal{F}_n\}$, to which we associate a class $\mathcal{M}_t$ of martingale estimating functions defined by

$$\mathcal{M}_t = \left\{ G_n(t) : G_n(t) = \sum_{j=1}^n w_j h_j(t), \ w_j = w_j(Y_1, \ldots, Y_{j-1}, \theta) \right\}. \tag{4}$$

Estimators of $\theta$ can then be found by solving the estimating equations $G_n(t) = 0$.

If we assume that the zero-mean martingale estimating functions $G_n(t)$ defined by (4) are square integrable and differentiable a.s. with respect to $\theta$ for each $t \in T$, then the special class of transform-based martingale estimating functions (4) fits into the general quasi-likelihood framework given by Godambe and Heyde [15]. In this framework, which incorporates essential ideas from the methods of least squares and maximum likelihood, a basic martingale $\{H_n = \sum_{j=1}^n h_j, \mathcal{F}_n\}$, with $h_j = h_j(Y_1, \ldots, Y_j, \theta)$ and $E(h_j|\mathcal{F}_{j-1}) = 0$ a.s., can be chosen in a variety of ways that give rise to different classes of martingale estimating functions as alternatives to the score function. In particular, since any discrete time process $\{Y_n, \mathcal{F}_n\}$ has the semimartingale representation $\sum_{j=1}^n Y_j = \sum_{j=1}^n E(Y_j|\mathcal{F}_{j-1}) + \sum_{j=1}^n h_j$, with $h_j = Y_j - E(Y_j|\mathcal{F}_{j-1})$, a class of martingale estimating functions may be based on the martingale $\sum_{j=1}^n h_j$.

In the present context, where it is assumed that a conditional transform $c_j(t) = E(g_t(Y_j)|\mathcal{F}_{j-1})$ can be readily obtained (as in the examples in Section 6), the martingale difference $h_j(t)$ defined in (3) leads to the semimartingale representation of the transformed process $g_t(Y_j)$, that is,

$$\sum_{j=1}^n g_t(Y_j) = \sum_{j=1}^n E(g_t(Y_j)|\mathcal{F}_{j-1}) + H_n(t), \qquad t \in T. \tag{5}$$

A strong law of large numbers for martingales will entail $H_n(t)/n \to 0$ a.s., for every $t \in T$. This asymptotic equivalence of $\sum_{j=1}^n g_t(Y_j)$ and $\sum_{j=1}^n E(g_t(Y_j)|$ $\mathcal{F}_{j-1}$) together with the one-to-one correspondence between the density $f(Y_j|Y_1, \ldots, Y_{j-1})$ and $E(g_t(Y_j)|\mathcal{F}_{j-1})$, for important kernels, supports using $\{H_n(t), \mathcal{F}_n\}$ as the basic martingale to generate the class $\mathcal{M}_t$ of estimating functions (4).

2.2. *Optimality considerations.* The general theory of quasi-likelihood furnishes an optimal estimating function within $\mathcal{M}_t$. Accordingly, the estimating function $G_n^*(t) \in \mathcal{M}_t$ given by

$$G_n^*(t) = \sum_{j=1}^n w_j^* h_j(t),$$

with

$$w_j^* = \frac{E(\frac{\partial}{\partial \theta} h_j(t)|\mathcal{F}_{j-1})}{E(h_j^2(t)|\mathcal{F}_{j-1})},$$

satisfies the small sample optimality criterion ($O_F$-optimality) of maximizing, for all $\theta$,

(6) $$\frac{[E(\frac{\partial}{\partial \theta} G_n(t))]^2}{E(G_n^2(t))},$$

and the asymptotic optimality criterion ($O_A$-optimality) of maximizing, a.s., for all $\theta$ and all $n \geq 1$,

(7) $$\frac{[\sum_{j=1}^n w_j E(\frac{\partial}{\partial \theta} h_j(t)|\mathcal{F}_{j-1})]^2}{\sum_{j=1}^n E((w_j h_j(t))^2|\mathcal{F}_{j-1})}.$$

The estimating function $G_n^*(t)$ is a quasi-score estimating function, and an estimator of $\theta$ obtained from $G_n^*(t) = 0$ is a quasi-likelihood estimator. A comprehensive explanation of quasi-likelihood concepts is available in Godambe and Heyde [15] and Heyde [19]. The quantity in (7), denoted by $I_{G_n(t)}$, is the martingale information in $G_n(t)$. Its maximum value, at $G_n(t) = G_n^*(t)$, is given by

$$I_{G_n^*(t)} = \sum_{j=1}^n E((w_j^* h_j(t))^2|\mathcal{F}_{j-1}).$$

$I_{G_n^*(t)}$ occurs as a scale variable in the asymptotic distribution of the quasi-likelihood estimator of $\theta$. For the score function $S_n$, $I_{S_n}$ is the conditional Fisher information $\sum_{j=1}^n E(s_j^2|\mathcal{F}_{j-1})$. Note that for $G_n^*(t)$, the quantity (6) is equal to $E(I_{G_n^*(t)})$. Explicit forms of $G_n^*(t)$ and $I_{G_n^*(t)}$ in terms of the kernel function may be obtained in view of $E(\frac{\partial}{\partial \theta} h_j(t)|\mathcal{F}_{j-1}) = -\frac{\partial}{\partial \theta} E(g_t(Y_j)|\mathcal{F}_{j-1})$ and $E(h_j^2(t)|\mathcal{F}_{j-1}) = \mathrm{Var}(g_t(Y_j)|\mathcal{F}_{j-1})$.

The choice of $h_j(t)$ as the martingale difference generates the family of classes $\mathcal{M} = \{\mathcal{M}_t, t \in T\}$ with a corresponding family $\mathcal{G}_n^* = \{G_n^*(t), t \in T\}$ of quasi-score functions and a family $\{\theta_t^*, t \in T\}$ of quasi-score estimators. We define a measure of conditional efficiency, $\text{eff}_c(G_n^*(t), S_n)$, for the quasi-score $G_n^*(t) \in \mathcal{M}_t$ relative to the score function $S_n$ as the ratio $I_{G_n^*(t)}/I_{S_n}$. Also, we define the conditional efficiency, $\text{eff}_c(G_n^*(t_1), G_n^*(t_2))$, of $G_n^*(t_1) \in \mathcal{M}_{t_1}$ relative to $G_n^*(t_2) \in \mathcal{M}_{t_2}$ as the ratio $I_{G_n^*(t_1)}/I_{G_n^*(t_2)}$. The efficiency of a quasi-score function $G_n^*$ may be defined alternatively in terms of the information quantity associated with $O_F$-optimality; see McLeish [26] and Merkouris [28]. However, because the quasi-score function $G_n^*$ and its information $I_{G_n^*}$ are expressed in terms of conditional functionals only, it is easier to compute $I_{G_n^*}$ than its unconditional counterpart $E(I_{G_n^*})$, which may not even exist (e.g., in models involving stable distributions). A justification for a nonasymptotic use of $I_{G_n^*}$ in measuring efficiency is provided in the next section.

We can now choose the most efficient quasi-score function in $\mathcal{G}_n^*$ by maximizing $I_{G_n^*(t)}$ with respect to $t \in T$. Thus, we will have

$$I_{G_n^*(t)} \leq I_{G_n^*(t^*)} \qquad \text{a.s., for some } t^* \in T,$$

and hence,

$$\text{eff}_c(G_n^*(t), S_n) \leq \text{eff}_c(G_n^*(t^*), S_n).$$

In general, the resulting estimator $\theta_{t^*}^*$ will be adaptive, in the sense that the value of $t^*$ will be determined by the sample. In the usual case where $I_{G_n^*(t)}$ depends on the parameter $\theta$, we may replace $\theta$ in $I_{G_n^*(t)}$ by an initial estimate and then proceed with the maximization; see Section 4.

The information contained in the sample $\{Y_j, 1 \leq j \leq n\}$, and carried by $\sum_{j=1}^n g_t(Y_j)$ for a specified kernel, is spread throughout the range of values of $t$. The above procedure of choosing the most efficient member of $\mathcal{G}_n^*$ aims to minimize the loss of information resulting from choosing any particular value for $t$. We may well then use more points from the set $T$ and extract the maximum information possible by judiciously choosing their values. This leads to the consideration of combining (in the sense of Heyde [18]) an arbitrary number of distinct transform martingale estimating functions. A distinctive advantage of the transform method is the ready capacity to form combinations of the form $\sum_{j=1}^n \sum_{l=1}^k w_{jl} h_j(t_l) = \sum_{j=1}^n \mathbf{w}_j \mathbf{h}_j(\mathbf{t})$ [in obvious notation for the vectors $\mathbf{w}_j$ and $\mathbf{h}_j(\mathbf{t})$] by using an arbitrary number of points $t_1, \ldots, t_k \in T$, and thereby producing more efficient transform quasi-score functions. A procedure for constructing transform-based composite martingale estimating functions that are statistically and computationally efficient is described in the next section.

Transform martingale estimating functions for a $p$-vector parameter $\boldsymbol{\theta}$ are $p$-dimensional and, for a $k$-vector $\mathbf{h}_j(\mathbf{t})$, have the general composite form $\mathbf{G}_n(\mathbf{t}) = \sum_{j=1}^{n} \mathbf{w}_j \mathbf{h}_j(\mathbf{t})$, in which the $\mathbf{w}_j$'s are $p \times k$ weighting matrices depending on $Y_1, \ldots, Y_{j-1}$ and $\boldsymbol{\theta}$. Suppressing $\mathbf{t}$ at the moment, the optimal $\mathbf{G}_n$ is $\mathbf{G}_n^* = \sum_{j=1}^{n} \mathbf{w}_j^* \mathbf{h}_j$, where $\mathbf{w}_j^* = (E(\dot{\mathbf{h}}_j | \mathcal{F}_{j-1}))'(E(\mathbf{h}_j \mathbf{h}_j' | \mathcal{F}_{j-1}))^{-1}$, $\dot{\mathbf{h}}_j = \{E(\frac{\partial}{\partial \theta_i} h_{jl})\}$ and prime denotes transpose. $\mathbf{G}_n^*$ maximizes, in the partial order of nonnegative definite matrices, the information matrix $\mathbf{I}_{\mathbf{G}_n} = \bar{\mathbf{G}}_n' \langle \mathbf{G} \rangle_n^{-1} \bar{\mathbf{G}}_n$, where $\langle \mathbf{G} \rangle_n = \sum_{j=1}^{n} E[(\mathbf{w}_j \mathbf{h}_j)(\mathbf{w}_j \mathbf{h}_j)' | \mathcal{F}_{j-1}]$ and $\bar{\mathbf{G}}_n = \sum_{j=1}^{n} \mathbf{w}_j E(\dot{\mathbf{h}}_j | \mathcal{F}_{j-1})$. Both $\langle \mathbf{G} \rangle_n$ and $\bar{\mathbf{G}}_n$ are assumed to be, a.s., nonsingular for each $n \geq 1$. For $\mathbf{G}_n^*$, it holds that $\mathbf{I}_{\mathbf{G}_n^*} = \langle \mathbf{G}^* \rangle_n = \sum_{j=1}^{n} E[(\mathbf{w}_j^* \mathbf{h}_j)(\mathbf{w}_j^* \mathbf{h}_j)' | \mathcal{F}_{j-1}]$. In analogy with a definition of efficiency, dating back to McLeish [26], which uses as measure of the size of the $O_F$ information matrix its determinant, we may define the conditional efficiency, $\mathrm{eff}_c(\mathbf{G}_n^*, \mathbf{S}_n)$, of $\mathbf{G}_n^*$ as the ratio $|\mathbf{I}_{\mathbf{G}_n^*}| / |\mathbf{I}_{\mathbf{S}_n}|$; elaboration on measures of efficiency based on the martingale information matrix can be found in Merkouris [28].

The procedure extends readily to multivariate observations. For $r$-dimensional random variables $\mathbf{Y}_j = (Y_{j1}, \ldots, Y_{jr})$, $1 \leq j \leq n$, the kernels are multivariate, with an $r$-dimensional index set, and are constructed as products of univariate kernels, that is,

$$(8) \qquad g_{\mathbf{t}}(\mathbf{Y}_j) = \prod_{i=1}^{r} g_{t_i}(Y_{ji}), \qquad \mathbf{t} = (t_1, \ldots, t_r) \in T^r.$$

2.3. *A link with a Fourier method for i.i.d. variables.* An interesting link of the proposed method of estimation with an existing transform-based method for i.i.d. variables is established as follows. Using the more suggestive notation $s_j = s(y_j; \mathcal{F}_{j-1})$, the score function $S_n = \sum_{j=1}^{n} s_j$ can be expressed as

$$(9) \qquad \sum_{j=1}^{n} s(y_j; \mathcal{F}_{j-1}) = \sum_{j=1}^{n} \int \omega_j(t) \exp(ity_j) \, dt,$$

where

$$\omega_j(t) = \frac{1}{2\pi} \int s(y_j; \mathcal{F}_{j-1}) \exp(-ity_j) \, dy_j, \qquad 1 \leq j \leq n,$$

is the inverse Fourier transform of $s(y_j; \mathcal{F}_{j-1})$. When the form of $\omega_j(t)$ is no more tractable than that of $s(Y_j; \mathcal{F}_{j-1})$, the integral in (9) may be approximated arbitrarily closely by a step function, say, $\sum_{l=1}^{k} w_j(t_l) \exp(it_l Y_j)$, the coefficients $w_{jl} = w_j(t_l)$ being functions of $Y_1, \ldots, Y_{j-1}$ and $\theta$ as well as $t$. Then the score function can be written as

$$
\begin{aligned}
(10) \qquad & \sum_{j=1}^{n} \int s(y; \mathcal{F}_{j-1}) \, d\hat{F}_j(y) \\
& = \sum_{j=1}^{n} \int \sum_{l=1}^{k} w_{jl} \exp(it_l y) \, d[\hat{F}_j(y) - F_j(y | \mathcal{F}_{j-1})],
\end{aligned}
$$

noticing that the term involving $F_j(y|\mathcal{F}_{j-1})$ is identically 0. In view of the transforms (1) and (2), the right-hand side of (10) can be written in the form

$$\sum_{j=1}^{n}\sum_{l=1}^{k}w_{jl}[\exp(it_lY_j)-E(\exp(it_lY_j)|\mathcal{F}_{j-1})]. \tag{11}$$

In the case of i.i.d. random variables, this approach leads to the generalized method of moments of Feuerverger and McDunnough [12, 13] based on the kernel $g_t(Y_j)=\exp(itY_j)$ and appropriate weights $w_{jl}$. A linear approximation of the score function for general classes of kernels was considered in the i.i.d. case by Feuerverger and McDunnough [14], and more extensively by Brant [4]. In the stochastic process context this leads to the generalization of (11),

$$
\begin{aligned}
\sum_{j=1}^{n}\sum_{l=1}^{k}&w_{jl}[g_{t_l}(Y_j)-E(g_{t_l}(Y_j)|\mathcal{F}_{j-1})]\\
&=\sum_{j=1}^{n}\sum_{l=1}^{k}w_{jl}h_j(t_l)=\sum_{l=1}^{k}\sum_{j=1}^{n}w_{jl}h_j(t_l),
\end{aligned}
\tag{12}
$$

which can be viewed as a combination of $k$ estimating functions of the form put forward earlier in this section. Thus, Fourier transform methods of estimation and more general transform-based linear approximations of the score function can be incorporated into the general quasi-likelihood theory.

**3. Combinations of transform martingale estimating functions.** In this section we develop a method of optimally combining transform-based martingales into quasi-score functions of nondecreasing conditional efficiency (as defined in Section 2.2). The construction of such optimal composite estimating functions, which for the main transforms can attain arbitrarily high efficiency, is founded on an alternative formulation of the optimality of a martingale estimating function based on the concept of orthogonal projection.

3.1. *Optimality and orthogonal projection.* Consider first the space $L^2 = L^2(\Omega, \mathcal{F}, P_{\boldsymbol{\theta}})$ of (equivalence classes of) random variables on $(\Omega, \mathcal{F}, P_{\boldsymbol{\theta}})$ which are square integrable (i.e., with finite second moment). Endowed with inner product $(X, Y) = E(XY)$ and norm $\|X\| = (X, X)^{1/2}$, the space $L^2$ is a Hilbert space. Let $\mathcal{A}$ be a closed subspace of $L^2$. For $X \in L^2$, let $E^*(X|\mathcal{A})$ denote the unique element in $\mathcal{A}$ such that

$$\|X - E^*(X|\mathcal{A})\|^2 = \inf_{Z \in \mathcal{A}}\|X - Z\|^2 = \inf_{Z \in \mathcal{A}}E[(X - Z)^2],$$

that is, $E^*(X|\mathcal{A})$ is the orthogonal projection of $X$ on $\mathcal{A}$.

Next consider the class $\mathcal{M}$ of general martingale estimating functions of the form $G_n = \sum_{j=1}^{n} \mathbf{w}_j \mathbf{h}_j$, generated by a basic martingale $\{\mathbf{H}_n = \sum_{j=1}^{n} \mathbf{h}_j, \mathcal{F}_n\}$, with $k$-dimensional martingale differences $\mathbf{h}_j = (h_{j1}, \ldots, h_{jk})'$ and with $k$-vector coefficients $\mathbf{w}_j = (w_{j1}, \ldots, w_{jk})$ that are $\mathcal{F}_{j-1}$-measurable functions depending on a scalar parameter $\theta$. The class $\mathcal{M}$ is a linear subspace of functions spanned by the $\mathbf{h}_j$'s. Henceforth, we will refer to the class $\mathcal{M}$ as a linear space, or simply a space.

A projection representation of the optimal martingale estimating function in the space $\mathcal{M}$ that appeared first in Merkouris [27] is formalized here in the following lemma. For a Hilbert space approach to general estimating functions based on the notion of $E$-sufficiency, see Small and McLeish [34].

LEMMA 1. *The martingale quasi-score function $G_n^* = \sum_{j=1}^{n} \mathbf{w}_j^* \mathbf{h}_j$, with $\mathbf{w}_j^* = (E(\dot{\mathbf{h}}_j|\mathcal{F}_{j-1}))'(E(\mathbf{h}_j\mathbf{h}_j'|\mathcal{F}_{j-1}))^{-1}$, is the orthogonal projection of the score function $S_n$ onto the space $\mathcal{M}$.*

PROOF. Let $\mathcal{M}_j$ denote the subspace of functions of the form $g_j = \mathbf{w}_j \mathbf{h}_j$, $1 \leq j \leq n$. Since the functions $\mathbf{h}_j$ are orthogonal, that is, $E(\mathbf{h}_i\mathbf{h}_j') = \mathbf{0}$ for $i \neq j$, the space $\mathcal{M}$ is the direct sum of the subspaces $\mathcal{M}_j$, that is,

$$(13) \qquad\qquad \mathcal{M} = \mathcal{M}_1 \oplus \cdots \oplus \mathcal{M}_n.$$

Now consider the Hilbert spaces $L^2(\Omega, \mathcal{F}_j, P_\theta^j)$, $1 \leq j \leq n$, where $P_\theta^j$ is the probability measure restricted to $\mathcal{F}_j$. Furthermore, consider the subspaces $B_j \subset L^2(\Omega, \mathcal{F}_j, P_\theta^j)$ of all measurable square integrable functions of $\mathbf{Y}_1, \ldots, \mathbf{Y}_j$ with conditional mean zero and differentiable with respect to $\theta$. Note here that the $j$th term of the score function, $s_j$, belongs to $B_j$ and the assumptions for $\mathcal{M}_j$ imply that $\mathcal{M}_j \subset B_j$.

Observe that since $E(\dot{\mathbf{h}}_j|\mathcal{F}_{j-1}) = E(s_j\mathbf{h}_j'|\mathcal{F}_{j-1})$ a.s., under a mild regularity condition $\mathbf{w}_j^*$ is the coefficient of $\dot{\mathbf{h}}_j$ in the orthogonal projection $g_j^* = \mathbf{w}_j^*\mathbf{h}_j$ of $s_j \in B_j$ onto the subspace $\mathcal{M}_j \subset B_j$. Specifically, denoting by $\|\cdot\|_{\mathcal{F}_{j-1}}$ the norm induced by $P_\theta^j$ conditional on $\mathcal{F}_{j-1}$, the element $g_j^* \in \mathcal{M}_j$ is such that

$$(14) \qquad\qquad \|s_j - g_j^*\|_{\mathcal{F}_{j-1}}^2 = \inf_{g_j \in \mathcal{M}_j} \|s_j - g_j\|_{\mathcal{F}_{j-1}}^2 \qquad \text{a.s.,}$$

for all $\theta \in \Theta$. The $\mathcal{F}_{j-1}$-measurable $\mathbf{w}_j^*$ arises then as the solution of the system of $k$ equations that express the orthogonality condition $E[(s_j - \mathbf{w}_j\mathbf{h}_j)\mathbf{h}_j'|\mathcal{F}_{j-1}] = \mathbf{0}$. By orthogonality, we also obtain $\|s_j\|_{\mathcal{F}_{j-1}}^2 = \|g_j^*\|_{\mathcal{F}_{j-1}}^2 + \|s_j - g_j^*\|_{\mathcal{F}_{j-1}}^2$ and, by passage to the sum, the decomposition of the conditional Fisher information $I_{S_n} = \sum_j \|s_j\|_{\mathcal{F}_{j-1}}^2$ into the (observed) information of the quasi-score function $I_{G_n^*} = \sum_j \|g_j^*\|_{\mathcal{F}_{j-1}}^2 = \sum_{j=1}^n E[(\mathbf{w}_j^* \mathbf{h}_j)(\mathbf{w}_j^* \mathbf{h}_j)'|\mathcal{F}_{j-1}]$ and the minimized sum of residuals $\sum_j \|s_j - g_j^*\|_{\mathcal{F}_{j-1}}^2$.

Now, since $s_j$ and $g_j$ are elements of $L^2(\Omega, \mathcal{F}_j, P_\theta^j)$, we can take the expectation of $\|s_j - g_j\|_{\mathcal{F}_{j-1}}^2$ to obtain

$$\|s_j - g_j^*\|^2 = \inf_{g_j \in \mathcal{M}_j} \|s_j - g_j\|^2$$

from (14), and so we can formally write $g_j^* = E^*(s_j|\mathcal{M}_j)$. Moreover, using the linearity property of the projection operator, the decomposition (13) and the fact that the martingale differences $s_i$ and $\mathbf{h}_j$, $i \neq j$, are mutually orthogonal, we have

$$E^*(S_n|\mathcal{M}) = \sum_{j=1}^n E^*(s_j|\mathcal{M}) = \sum_{j=1}^n E^*(s_j|\mathcal{M}_j) = \sum_{j=1}^n g_j^* = G_n^*.$$

Therefore, the quasi-score function $G_n^*$ is the unique orthogonal projection of the score function $S_n$ onto $\mathcal{M}$. $\square$

A few remarks are in order.

REMARK 1. As expectation conditional on $\mathcal{F}_{j-1}$, the inner product that induces the norm $\|\cdot\|_{\mathcal{F}_{j-1}}$ used in (14) is an element of the space of all $\mathcal{F}_{j-1}$-measurable functions, say, $\mathcal{M}_{\mathcal{F}_{j-1}}$. Such an inner product, with the associated space of "scalars" being $\mathcal{M}_{\mathcal{F}_{j-1}}$, is well defined, its defining properties holding a.s. For a comparable inner product used in a similar context, see Murphy and Li [29]. It is important to note here that if the $g_j^*$s are only conditionally square integrable (see Example 1 in Section 6), then the result of the lemma holds, but restricted to term-wise projection of the score function onto the $\mathcal{M}_j$'s using the norm $\|\cdot\|_{\mathcal{F}_{j-1}}$. In this case the $O_F$ information quantity $E(I_{G_n^*})$ does not exist.

REMARK 2. For a $p$-vector parameter $\boldsymbol{\theta}$, each component $s_{ji}$ of $\mathbf{s}_j$, $i = 1, \ldots, p$, is approximated by its projection $g_{ji}^* = \mathbf{w}_{ji}^* \mathbf{h}_j$ onto the same space $\mathcal{M}_j$. More compactly, we write $\mathbf{g}_j^* = \mathbf{w}_j^* \mathbf{h}_j$, where the $p \times k$ matrix $\mathbf{w}_j^* = (\mathbf{w}_{j1}^*, \ldots, \mathbf{w}_{jp}^*)$ is given by $\mathbf{w}_j^* = (E(\dot{\mathbf{h}}_j|\mathcal{F}_{j-1}))'(E(\mathbf{h}_j \mathbf{h}_j'|\mathcal{F}_{j-1}))^{-1}$. Both $\mathbf{g}_j^*$ and $\mathbf{s}_j = (s_{j1}, \ldots, s_{jp})'$ are elements of the set $L_2^p$ of random $p$-vectors with all components in $L^2(\Omega, \mathcal{F}_j, P_{\boldsymbol{\theta}}^j)$, which is a Hilbert space when the inner product is defined to be $(\mathbf{X}, \mathbf{Y}) = \operatorname{tr} E(\mathbf{X}\mathbf{Y}')$ for all $\mathbf{X}, \mathbf{Y} \in L_2^p$. In this Hilbert space, the vector $\mathbf{g}_j^*$ can be characterized as the projection of $\mathbf{s}_j$ onto $\mathcal{M}_j^p$, the $p$-fold Cartesian product of $\mathcal{M}_j$ with itself, and $\mathbf{G}_n^*$ as the projection of $\mathbf{S}_n$ onto $\mathcal{M}^p = \mathcal{M}_1^p \oplus \cdots \oplus \mathcal{M}_n^p$.

REMARK 3.    Enlarging the space of martingale estimating functions increases the information of the corresponding martingale quasi-score function. For an increasing sequence of spaces $\{\mathcal{M}_k\}$ of martingale estimating functions with corresponding sequence of quasi-score functions $\{G^*_{n,k}\}$, we have $G^*_{n,k+1} = E^*(S_n|\mathcal{M}_{k+1})$ and, by a well-known projection property,

$$E^*(G^*_{n,k+1}|\mathcal{M}_k) = E^*(E^*(S_n|\mathcal{M}_{k+1})|\mathcal{M}_k) = E^*(S_n|\mathcal{M}_k) = G^*_{n,k},$$

that is, $G^*_{n,k}$ is the orthogonal projection of $G^*_{n,k+1}$ onto $\mathcal{M}_k$, and $\{G^*_{n,k}, \mathcal{M}_k\}$ is a "projection martingale." It follows easily that $\|S_n - G^*_{n,k+1}\|^2 \leq \|S_n - G^*_{n,k}\|^2$. In view of (14), this also entails $I_{G^*_{n,k}} \leq I_{G^*_{n,k+1}}$ a.s.

REMARK 4.    Combining martingale estimating functions essentially increases the dimensionality of $\mathbf{h}_j$, for each $j$, and hence, of the space $\mathcal{M}_j$ spanned by its components. Thus, according to the previous remark, a closer approximation of the score function can be achieved. In the furthermost extension of the combination procedure, aiming for full efficiency, we consider the closure $\bar{\mathcal{M}}_j$ of the infinite-dimensional space $\mathcal{M}_j$ spanned by a countable set $\{h_{j1}, h_{j2}, \ldots\}$ of martingale differences. The limits associated with $\bar{\mathcal{M}}_j$ are in the norm $\|\cdot\|_{\mathcal{F}_{j-1}}$. If the set $\{h_{j1}, h_{j2}, \ldots\}$ is complete in $B_j \subset L^2(\Omega, \mathcal{F}_j, P^j_\theta)$, that is, if $\bar{\mathcal{M}}_j = B_j$, for all $1 \leq j \leq n$, then

$$g^*_j = E^*(s_j|\bar{\mathcal{M}}_j) = E^*(s_j|B_j) = s_j,$$

and hence $G^*_n = S_n$. This possibility is explored next using transform martingale estimating functions.

3.2. *Optimal combinations of transform martingales.*  We assume, for clarity, that $\theta$ is scalar and recall that the most efficient quasi-score function $G^*_n(t^*)$ in the family $\mathcal{G}^*_n = \{G^*_n(t), t \in T\}$ of quasi-score functions corresponding to the family of distinct spaces $\mathcal{M} = \{\mathcal{M}_t, t \in T\}$ is obtained by maximizing $I_{G^*_n(t)}$ with respect to $t \in T$. The construction of composite martingale estimating functions with the use of more points $t \in T$ involves increasing the dimension of the basic transform martingale and, hence, the dimension of the associated space of transform martingale estimating functions. The proposed procedure generates in a stepwise manner an increasing sequence of spaces of transform martingale estimating functions by retaining the optimal points $t \in T$ determined in the preceding steps. This facilitates the comparison of composite martingale quasi-scores from different spaces and ensures increasing efficiency with increasing number of points $t \in T$.

Proceeding formally, we write $t^*_l$, $1 \leq l \leq k-1$, $k \geq 2$, for the optimal points in $T$ determined in the first $k-1$ steps. At the $k$th step of the approximation, with a new point $t_k \neq t^*_l$, we will have the $k \times 1$ vector martingale $\{H_n(t^*_1, \ldots, t^*_{k-1}, t_k) = (\sum_{j=1}^n h_j(t^*_1), \ldots, \sum_{j=1}^n h_j(t^*_{k-1}), \sum_{j=1}^n h_j(t_k))', \mathcal{F}_n\}$ with associated space $\mathcal{M}_{t_1^*,\dots,t_{k-1}^*,t_k}$ of composite martingale estimating functions of the form

$$G_n(t_1^*,\dots,t_{k-1}^*,t_k) = \sum_{j=1}^n \mathbf{w}_j \mathbf{h}_j = \sum_{j=1}^n \left( \sum_{l=1}^{k-1} w_{jl} h_j(t_l^*) + w_{jk} h_j(t_k) \right),$$

with $t_l^*$ fixed. Then the quasi-score function within $\mathcal{M}_{t_1^*,\dots,t_{k-1}^*,t_k}$ is given by

$$G_n^*(t_1^*,\dots,t_{k-1}^*,t_k) = \sum_{j=1}^n \mathbf{w}_j^* \mathbf{h}_j,$$

with $\mathbf{w}_j^* = (E(\dot{\mathbf{h}}_j | \mathcal{F}_{j-1}))'(E(\mathbf{h}_j \mathbf{h}_j' | \mathcal{F}_{j-1}))^{-1}$, and the information in $G_n^*(t_1^*,\dots,t_{k-1}^*,t_k)$ is

$$I_{G_n^*(t_1^*,\dots,t_{k-1}^*,t_k)} = \sum_{j=1}^n (E(\dot{\mathbf{h}}_j | \mathcal{F}_{j-1}))'(E(\mathbf{h}_j \mathbf{h}_j' | \mathcal{F}_{j-1}))^{-1}(E(\dot{\mathbf{h}}_j | \mathcal{F}_{j-1})).$$

Clearly, $\mathcal{M}_{t_1^*,\dots,t_{k-1}^*} \subseteq \mathcal{M}_{t_1^*,\dots,t_{k-1}^*,t_k}$ for any $t_k \in T$, and since $G_n^*(t_1^*,\dots,t_{k-1}^*,t_k)$ is the projection of $S_n$ onto $\mathcal{M}_{t_1^*,\dots,t_{k-1}^*,t_k}$, we have $I_{G_n^*(t_1^*,\dots,t_{k-1}^*)} \leq I_{G_n^*(t_1^*,\dots,t_{k-1}^*,t_k)}$ a.s. for any $t_k \in T$ (by Remark 3, Section 3.1). From the family $\{G_n^*(t_1^*,\dots,t_{k-1}^*,t_k), t_k \in T\}$ we can now choose the most efficient quasi-score function by maximizing $I_{G_n^*(t_1^*,\dots,t_{k-1}^*,t_k)}$ with respect to $t_k \in T$. Thus, for some $t_k^* \in T$, we will have

$$(15) \qquad G_n^*(t_1^*,\dots,t_k^*) = E^*(S_n | \mathcal{M}_{t_1^*,\dots,t_k^*}) = \sum_{j=1}^n \sum_{l=1}^k w_j^*(t_l) h_j(t_l^*),$$

satisfying

$$(16) \qquad I_{G_n^*(t_1^*,\dots,t_{k-1}^*)} \leq I_{G_n^*(t_1^*,\dots,t_{k-1}^*,t_k)} \leq I_{G_n^*(t_1^*,\dots,t_{k-1}^*,t_k^*)} \qquad \text{a.s.,}$$

and hence

$$\text{eff}_c(G_n^*(t_1^*,\dots,t_{k-1}^*), S_n) \leq \text{eff}_c(G_n^*(t_1^*,\dots,t_k^*), S_n) \qquad \text{a.s.}$$

Moreover,

$$(17) \qquad \text{eff}_c(G_n^*(t_1^*,\dots,t_k^*), S_n) \leq \text{eff}_c(G_n^*(t_1^*,\dots,t_k^*), G_n^*(t_1^*,\dots,t_{k+1}^*)) \qquad \text{a.s.}$$

It is worth noting that the inequality (16) provides a nondecreasing sequence of lower bounds for the score information $I_{S_n}$, while the inequality (17) places an upper bound to the unknown efficiency of $G_n^*(t_1^*,\dots,t_k^*)$. The procedure of adding components to $G_n^*(t_1^*,\dots,t_k^*)$ may be terminated when $I_{G_n^*(t_1^*,\dots,t_{k-1}^*,t_k^*)}$ ceases to increase substantially with $k$ or when computations become prohibitive.

The composite quasi-score function in (15) provides a finite approximation to the score function $S_n$ in $L^2(\Omega, \mathcal{F}, P_\theta)$. We may consider next an arbitrarily close approximation of $S_n$ by projection onto infinite-dimensional spaces of transform martingale estimating functions. Thus, for each $j$, let $\{h_j(t_1), h_j(t_2), \ldots\}$ be a countable set formed by evaluating the martingale difference $h_j(\cdot)$ at distinct points $t_1, t_2, \ldots$ in $T$. Let $\bar{\mathcal{M}}_{j;t_1,t_2,\ldots}$ be the closed linear subspace spanned by this set. $\bar{\mathcal{M}}_{j;t_1,t_2,\ldots}$ consists of those zero-mean elements of $L^2(\Omega, \mathcal{F}_j, P_\theta^j)$ which can be approximated in the norm $\|\cdot\|_{\mathcal{F}_{j-1}}$ by finite linear combinations of elements of $\{h_j(t_1), h_j(t_2), \ldots\}$ with $\mathcal{F}_{j-1}$-measurable coefficients. Note here that a set complete in $L^2(\Omega, \mathcal{F}_j, P_\theta^j)$ is characterized by the property that $L^2(\Omega, \mathcal{F}_j, P_\theta^j)$ contains no nonzero function orthogonal to all elements of the set (see Burrill [6], page 216). It follows then that the set $\{h_j(t_1), h_j(t_2), \ldots\}$ is complete a.s. in the subspace $B_j \subset L^2(\Omega, \mathcal{F}_j, P_\theta^j)$, that is, $\bar{\mathcal{M}}_{j;t_1,t_2,\ldots} = B_j$ a.s., if the set $\{g_{t_1}(Y_j), g_{t_2}(Y_j), \ldots\}$ is complete in $L^2(\Omega, \mathcal{F}_j, P_\theta^j)$ in the norm $\|\cdot\|$. This leads to the following theorem.

THEOREM 1.   *Suppose that for a kernel class of functions $\{g_t(Y), t \in T\}$ the countable set $\{g_{t_1}(Y_j), g_{t_2}(Y_j), \ldots\}$ is complete in $L^2(\Omega, \mathcal{F}_j, P_\theta^j)$ for all $j$. Then the score function $S_n$ can be approximated in $L^2(\Omega, \mathcal{F}, P_\theta)$ arbitrarily closely by the transform quasi-score function $G_n^*(t_1, \ldots, t_k)$ with $k$ sufficiently large.*

PROOF.   Since $s_j \in B_j$, the completeness of the set $\{g_{t_1}(Y_j), g_{t_2}(Y_j), \ldots\}$ in $L^2(\Omega, \mathcal{F}_j, P_\theta^j)$ implies $s_j \in \bar{\mathcal{M}}_{j;t_1,t_2,\ldots}$ a.s. Then $s_j$ is a limit point in $\bar{\mathcal{M}}_{j;t_1,t_2,\ldots}$ and, thus, for each $\varepsilon > 0$ there exists a finite combination $g_j(t_1, t_2, \ldots, t_K) = \sum_{l=1}^K w_{jl} h_j(t_l)$, an element of the subspace $\mathcal{M}_{j;t_1,t_2,\ldots,t_K}$, for which $\|s_j - g_j(t_1, t_2, \ldots, t_K)\|_{\mathcal{F}_{j-1}}^2 < \varepsilon$   a.s. But by (14),

$$\|s_j - g_j^*(t_1, t_2, \ldots, t_K)\|_{\mathcal{F}_{j-1}}^2 \leq \|s_j - g_j(t_1, t_2, \ldots, t_K)\|_{\mathcal{F}_{j-1}}^2 < \varepsilon, \qquad \text{a.s.,}$$

for the projection $g_j^*(t_1, t_2, \ldots, t_K) = \sum_{l=1}^K w_{jl}^* h_j(t_l)$ of $s_j$ onto $\mathcal{M}_{j;t_1,t_2,\ldots,t_K}$. Now, by Remark 3, Section 3.1,

$$\|s_j - g_j^*(t_1, t_2, \ldots, t_{K+1})\|_{\mathcal{F}_{j-1}}^2 \leq \|s_j - g_j^*(t_1, t_2, \ldots, t_K)\|_{\mathcal{F}_{j-1}}^2 \qquad \text{a.s.,}$$

and hence

$$\|s_j - g_j^*(t_1, t_2, \ldots, t_k)\|_{\mathcal{F}_{j-1}}^2 < \varepsilon \qquad \text{a.s., for all } k \geq K.$$

This implies

$$\|s_j - g_j^*(t_1, t_2, \ldots, t_k)\|_{\mathcal{F}_{j-1}}^2 \to 0 \qquad \text{a.s., as } k \to \infty,$$

and by the decomposition $\|s_j\|^2_{\mathcal{F}_{j-1}} = \|g^*_j(t_1, t_2, \ldots, t_k)\|^2_{\mathcal{F}_{j-1}} + \|s_j - g^*_j(t_1, t_2, \ldots, t_k)\|^2_{\mathcal{F}_{j-1}}$ it also implies

$$(18) \qquad \|g^*_j(t_1, t_2, \ldots, t_k)\|^2_{\mathcal{F}_{j-1}} \to \|s_j\|^2_{\mathcal{F}_{j-1}} \qquad \text{a.s.}$$

By the (unconditional) square integrability of the elements $s_j$ and $g^*_j(t_1, t_2, \ldots, t_k)$ of $L^2(\Omega, \mathcal{F}_j, P^j_\theta)$, we have

$$\|s_j - g^*_j(t_1, t_2, \ldots, t_k)\|^2 \to 0 \qquad \text{as } k \to \infty,$$

or, equivalently,

$$s_j = g^*_j(t_1, t_2, \ldots) = \sum_{l=1}^{\infty} w^*_{jl} h_j(t_l).$$

It follows that

$$S_n = G^*_n(t_1, t_2, \ldots) = \sum_{j=1}^{n} \sum_{l=1}^{\infty} w^*_{jl} h_j(t_l),$$

where $G^*_n(t_1, t_2, \ldots)$ is the quasi-score function within the space $\bar{\mathcal{M}}_{t_1, t_2, \ldots} \equiv \bar{\mathcal{M}}_{1; t_1, t_2, \ldots} \oplus \cdots \oplus \bar{\mathcal{M}}_{n; t_1, t_2, \ldots}$, that is, the unique orthogonal projection $E^*(S_n | \bar{\mathcal{M}}_{t_1, t_2, \ldots})$. Therefore, the score function $S_n$ can be approximated arbitrarily closely by the partial sums $G^*_n(t_1, \ldots, t_k) = \sum_{j=1}^{n} \sum_{l=1}^{k} w^*_{jl} h_j(t_l)$.

Moreover, by (18),

$$I_{S_n} = \lim_{k \to \infty} I_{G^*_n(t_1, t_2, \ldots, t_k)} \qquad \text{a.s.},$$

so that $\lim_{k \to \infty} \text{eff}_c(G^*_n(t_1, t_2, \ldots, t_k), S_n) = 1$ a.s.   $\square$

REMARK 5.    Completeness ensures full efficiency of the transform quasi-score function $G^*_n(t_1, t_2, \ldots)$, but is not necessary. It would suffice that $S_n \in \mathcal{M}_{t_1, \ldots, t_k}$ for some $k$ and $t_1, \ldots, t_k$ in $T$, or $S_n \in \bar{\mathcal{M}}_{t_1, t_2, \ldots}$. The completeness of kernel sets for the main transforms is shown in Feuerverger and McDunnough [14].

REMARK 6.    In the case that the probability measures $P^j_\theta$ are concentrated on a finite set of points, say, $\{a_1, \ldots, a_N\}$, the spaces $L^2(\Omega, \mathcal{F}_j, P^j_\theta)$ are $N$-dimensional, and then any finite set $\{g_{t_1}, \ldots, g_{t_N}\}$ of linearly independent kernel functions will be complete. Linear independence of the kernel functions requires that the vectors $(g_{t_l}(a_1), \ldots, g_{t_l}(a_N))'$, $1 \le l \le N$, be independent for some values of $t_1, \ldots, t_N$. Any such choice of values is optimal and yields $G^*_n(t_1, \ldots, t_N) \equiv S_n$.

REMARK 7. In general, full efficiency may be achieved with a finite number of points in $T$ for some terms of $S_n$, that is, $s_j = \sum_{l=1}^{k} w_{jl}^* h_j(t_l)$, if $s_j \in \mathcal{M}_{j;t_1,\ldots,t_k}$ for some $j$. In particular, $s_j = w_j^* h_j(t)$ may hold for some $t \in T$ when the conditional density $f(Y_j|Y_1,\ldots,Y_{j-1})$ belongs to the exponential family. This can occur, for example, in applications to aggregate Markov chains. In such situations projecting $s_j$ onto the larger space of functions of the form $g_j(t_1,t_2) = w_{j1}h_j(t_1) + w_{j2}h_j(t_2)$ results in $s_j = w_{j1}^* h_j(t_1) + 0 \cdot h_j(t_2)$ and $I_{g_j^*(t_1,t_2)} = I_{g_j^*(t_1)}$. It is not difficult to show that this is equivalent to

$$E(\dot{h}_j(t_1)|\mathcal{F}_{j-1})E(h_j(t_1),h_j(t_2)|\mathcal{F}_{j-1}) = E(\dot{h}_j(t_2)|\mathcal{F}_{j-1})E(h_j^2(t_1)|\mathcal{F}_{j-1}),$$

which can be written as

$$\frac{\partial}{\partial\theta}E(g_{t_1}(Y_j)|\mathcal{F}_{j-1})\operatorname{Cov}(g_{t_1}(Y_j),g_{t_2}(Y_j)|\mathcal{F}_{j-1})$$

$$= \frac{\partial}{\partial\theta}E(g_{t_2}(Y_j)|\mathcal{F}_{j-1})\operatorname{Var}(g_{t_1}(Y_j)|\mathcal{F}_{j-1}).$$

This condition is necessary for $f(Y_j|Y_1,\ldots,Y_{j-1})$ to belong to the exponential family with $g_{t_1}(Y_j)$ being the sufficient statistic for $\theta$. If, as in some situations, the sufficient statistic is $g_{t_1}(Y_j) = Y_j^{t_1}$, $t_1 = 1$, then

$$s_j = w_{j1}^* h_j(t_1) = -\frac{\frac{\partial}{\partial\theta}E(Y_j|\mathcal{F}_{j-1})}{\operatorname{Var}(Y_j|\mathcal{F}_{j-1})}[Y_j - E(Y_j|\mathcal{F}_{j-1})].$$

The multiparameter case can be treated in a similar manner. It is important to note that, since the sequence of the constructed spaces is increasing, a matrix version of the inequality (16) holds in the partial order of nonnegative matrices. This enables us to choose at the $k$th step of the procedure the optimal quasi-score $\mathbf{G}_n^*(t_1^*,\ldots,t_{k-1}^*,t_k)$ by maximizing the determinant of the matrix $\mathbf{I}_{\mathbf{G}_n^*(t_1^*,\ldots,t_{k-1}^*,t_k)}$ with respect to $t_k$. In the case of multivariate (say, $r$-dimensional) observations, the $k \times 1$ vector martingale difference $\mathbf{h}_j = (h_j(\mathbf{t}_1),\ldots,h_j(\mathbf{t}_k))'$ involves $k$ $r$-tuples in $T^r$. Apart from the computational difficulties associated with the choice of appropriate values of $\mathbf{t}_1,\ldots,\mathbf{t}_k$ for large $r$, the described methodology carries over to this case in a straightforward manner.

## 4. Kernel classes and computational issues.

The statistical and computational efficiency of the transform method depends on the associated kernel, the number $k$ of points $t_1,\ldots,t_k$, and the choice of the values of these points at which the kernel and its conditional expectation are evaluated.

Conditional transforms conveniently describing stochastic process models are, commonly, conditional versions of the characteristic function, the moment generating function, the Laplace transform and the probability generating function. The sequence of conditional moments can be derived from these transforms. As Brant [4] notes in the context of i.i.d. variables, the corresponding kernel classes are closed under multiplication, that is,

$$(19) \qquad g_t(Y)g_s(Y) = g_{v(t,s)}(Y) \qquad \text{for all } t, s \text{ in } T,$$

with a multiplication rule $v : T \times T \to T$ defined by (19) for any particular class. For instance, in the case of the moment generating function kernel $g_t(Y) = \exp(tY)$, we have $g_t(Y)g_s(Y) = g_{t+s}(Y)$, for all $t$, $s$ in $T$. The important practical consequence of the closure property of the kernel classes is that the transform $c_j(t) = E(g_t(Y_j)|\mathcal{F}_{j-1})$ yields the joint moment structure $c_j(v(t,s))$ of the kernel class, for $1 \le j \le n$. Thus, the matrix $E(\mathbf{h}_j \mathbf{h}_j' | \mathcal{F}_{j-1})$, with $\mathbf{h}_j = (h_j(t_1), \ldots, h_j(t_k))'$, can be readily obtained since its entries are of the form

$$E(h_j(t_i)h_j(t_l)|\mathcal{F}_{j-1}) = \text{Cov}(g_{t_i}(Y_j), g_{t_l}(Y_j)|\mathcal{F}_{j-1})$$
$$= c_j(v(t_i, t_l)) - c_j(t_i)c_j(t_l), \qquad i \ne l,$$

for any of the aforementioned kernels, except for the kernel $g_t(Y) = \exp(itY)$ of the characteristic function for which closure is under complex conjugation and $c_j(v(t_i, t_l)) - c_j(t_i) = c_j(t_i - t_l) - c_j(-t_i)$. Therefore, optimal combinations of transform martingale estimating functions can be readily constructed. The property of closure of the kernel classes under multiplication is preserved in multivariate kernels in an obvious way, in view of their defining property (8).

A special feature of the characteristic function with kernel class $\{g_t(Y) = \exp(itY), t \in \mathbb{R}\}$ merits attention. Notwithstanding the notational convenience associated with using the complex valued kernel, for computational purposes it is preferable to work with the class of real valued kernel vectors $\{(\cos(tY), \sin(tY)), t \in \mathbb{R}\}$. Thus, we may start with writing the basic martingale difference $h_j(t) = g_t(Y_j) - E(g_t(Y_j)|\mathcal{F}_{j-1})$ in the complex domain form

$$h_j(t) = \cos(tY_j) - E(\cos(tY_j)|\mathcal{F}_{j-1}) + i[\sin(tY_j) - E(\sin(tY_j)|\mathcal{F}_{j-1})].$$

Then, for the real valued vector $\mathbf{g}(Y_j) = (\cos(t_1 Y_j), \ldots, \cos(t_k Y_j), \sin(t_1 Y_j), \ldots, \sin(t_k Y_j))'$, the martingale difference vector is $\mathbf{h}_j = (\text{Re } h_j(t_1), \ldots, \text{Re } h_j(t_k), \text{Im } h_j(t_1), \ldots, \text{Im } h_j(t_k))'$. As in the i.i.d. case (e.g., Brant [4]), a convenient choice of kernel vector is $\mathbf{g}(Y_j) = (\cos(\tau Y_j), \ldots, \cos(k\tau Y_j), \sin(\tau Y_j), \ldots, \sin(k\tau Y_j))'$, for some $\tau$, which provides an approximation to the score function by a trigonometric polynomial.

For specified kernel and fixed $k$, the transform method requires the choice of the values of $\{t_l\}$. The most convenient approach is to choose the values a priori with a uniform spacing suitable to the particular problem. This has been tried with other transform methods in the case of i.i.d. variables and in connection mainly with the characteristic function; see, for example,

Epps and Pulley [9]. The efficiency of the procedure may vary with different uniform spacings. In the present context, an optimal uniform spacing could be sought by maximizing the martingale information with respect to the distance, say, $\tau$, between the points.

It may be noted parenthetically that the sequence of moments, with kernel class $\{g_t(Y) = Y^t,\ t = 0, 1, 2, \ldots\}$, has the distinct advantage over the other transforms that there is a natural choice of points $\{t_l\}$ as the first $k$ positive integers. The first $k$ should normally be chosen since the variability of the sample moments increases with their order. For $\mathbf{g}(Y_j) = (Y_j, Y_j^2, \ldots, Y_j^k)$, the matrix $E(\mathbf{h}_j \mathbf{h}_j' | \mathcal{F}_{j-1})$ involves moments up to order $2k$. The transform quasi-score function based on integer moments is a polynomial and, thus, its efficiency depends on how close the terms of the score function are to a polynomial. We will be denoting the polynomial quasi-score function based on the first $k$ moments by $G_n^*(1, \ldots, k)$.

In general, for small $k$ an optimal arbitrary spacing may result in sharp improvement in the efficiency. The values of $\{t_l\}$ can be chosen so that they maximize the martingale information. This can be taken as the "optimal choice" rule with respect to efficiency. However, as $k$ increases, so do the computational complexities of optimizing a $k$-dimensional surface and inverting $E(\mathbf{h}_j \mathbf{h}_j' | \mathcal{F}_{j-1})$. It is more convenient, following the procedure described in the previous section, to form the increasing sequence of spaces $\mathcal{M}_{t_1^*, \ldots, t_{k-1}^*} \subseteq \mathcal{M}_{t_1^*, \ldots, t_{k-1}^*, t_k}$, $k \geq 2$, holding the first $k-1$ points fixed at the values $t_1^*, \ldots, t_{k-1}^*$, and then to maximize $I_{G_n^*(t_1^*, \ldots, t_{k-1}^*, t_k)}$ with respect to $t_k$. Maximizing the information with respect to one point at a time may result in some loss of efficiency, which in many examples becomes negligible as $k$ grows larger.

The appropriate values of $\{t_l\}$ will generally depend on the unknown parameter. We propose a two-step approach in which we start with $I_{G_n^*(t)}$ evaluated at some preliminary estimate of the parameter (e.g., conditional least squares estimate), and at the $k$th stage of the approximation we evaluate $I_{G_n^*(t_1^*, \ldots, t_{k-1}^*, t_k)}$ at the value of the estimate obtained from $G_n^*(t_1^*, \ldots, t_{k-1}^*) = 0$. Then we use the value of the optimal point $t_k^*$ to solve the estimating equation $G_n^*(t_1^*, \ldots, t_{k-1}^*, t_k^*) = 0$. The solution is now the updated estimate in $I_{G_n^*(t_1^*, \ldots, t_{k-1}^*, t_k)}$, and a new value of the optimal point $t_k^*$ can be determined. The iteration continues until convergence to some value of the optimal point. This value of the optimal point is used in the estimating equation $G_n^*(t_1^*, \ldots, t_{k-1}^*, t_k^*) = 0$ to obtain the transform quasi-likelihood estimate and the value of the information quantity. This iterative scheme converges more rapidly as the number of points $\{t_l\}$, and hence, the efficiency of the estimating function, increases. In light of this, as more points $t$ are introduced the iteration may be stopped after the first step to ease the computations.

**5. Comparing transform estimating functions.** In this section we compare the polynomial quasi-score function, based on moments, with quasi-score functions based on other important transforms, namely, the characteristic function (CF), the moment generating function (MGF) and the probability generating function (PGF). This comparison is of theoretical and of methodological interest.

5.1. *Comparing $G_n^*(t_1, \ldots, t_k)$ with $G_n^*(1, \ldots, k)$.* A relationship of the CF, MGF and PGF quasi-score functions with the moment (M) quasi-score function is established by the following proposition.

PROPOSITION 1. *Assume that $E(Y_j^{2k}|\mathcal{F}_{j-1})$ exists. Then for both the CF transform, with the complex valued kernel $\exp(itY)$, and the MGF transform, the following relationship holds:*

$$\lim_{\max\{|t_l|\} \to 0} G_n^*(t_1, \ldots, t_k) = G_n^*(1, \ldots, k).$$

*For the PGF transform,*

$$\lim_{\max\{|t_l - 1|\} \to 0} G_n^*(t_1, \ldots, t_k) = G_n^*(1, \ldots, k).$$

PROOF. For $k = 1$, on the assumption that differentiation with respect to $t$ and expectation operations can be interchanged in

$$G_n^*(t) = -\sum_{j=1}^n \frac{\frac{\partial}{\partial\theta}E(g_t(Y_j)|\mathcal{F}_{j-1})}{\text{Var}(g_t(Y_j)|\mathcal{F}_{j-1})}[g_t(Y_j) - E(g_t(Y_j)|\mathcal{F}_{j-1})],$$

taking the limit and applying l'Hospital's rule (twice) yields the limit

$$G_n^*(1) = -\sum_{j=1}^n \frac{\frac{\partial}{\partial\theta}E(Y_j|\mathcal{F}_{j-1})}{\text{Var}(Y_j|\mathcal{F}_{j-1})}[Y_j - E(Y_j|\mathcal{F}_{j-1})].$$

Alternatively, a first-order Taylor series expansion of $g_t(Y_j)$ gives a representation of $G_n^*(t)$ in terms of $G_n^*(1)$ as

$$(20) \qquad G_n^*(t) = G_n^*(1) + t\sum_{j=1}^n K_j + o(t),$$

where

$$K_j = \left[\frac{\partial}{\partial\theta}E(Y_j|\mathcal{F}_{j-1})\,\text{Var}(Y_j|\mathcal{F}_{j-1})[Y_j^2 - E(Y_j^2|\mathcal{F}_{j-1})]\right.$$

$$\left. + \left[\frac{\partial}{\partial\theta}E(Y_j^2|\mathcal{F}_{j-1})\,\text{Var}(Y_j|\mathcal{F}_{j-1})\right.$$

$$-2\frac{\partial}{\partial\theta}E(Y_j|\mathcal{F}_{j-1})\operatorname{Cov}(Y_j,Y_j^2|\mathcal{F}_{j-1})\Big]$$

$$\times[Y_j-E(Y_j|\mathcal{F}_{j-1})]\Big]\Big/2[\operatorname{Var}(Y_j|\mathcal{F}_{j-1})]^2.$$

If $t$ is small and the $o(t)$ terms are neglected, then the quasi-score $G_n^*(t)$ appears as a perturbed version of $G_n^*(1)$. We have also

$$(21) \qquad I_{G_n^*(t)}=I_{G_n^*(1)}+t\sum_{j=1}^{n}L_j+o(t),$$

where

$$L_j=\frac{\frac{\partial}{\partial\theta}E(Y_j|\mathcal{F}_{j-1})}{\operatorname{Var}^2(Y_j|\mathcal{F}_{j-1})}\Big[\frac{\partial}{\partial\theta}E(Y_j^2|\mathcal{F}_{j-1})\operatorname{Var}(Y_j|\mathcal{F}_{j-1})$$

$$-\frac{\partial}{\partial\theta}E(Y_j|\mathcal{F}_{j-1})\operatorname{Cov}(Y_j,Y_j^2|\mathcal{F}_{j-1})\Big].$$

For $k\geq 2$, the results can be obtained similarly with the help of a symbolic computation package. □

Note that for the CF transform with the real valued kernel $(\cos(tY),\sin(tY))$ the limiting quasi-score is $G_n^*(1,\ldots,2k)$. If $m<k$ points $t$ tend to zero (or to one for the PGF transform), the limiting quasi-score function involves a kernel vector whose corresponding $m$ components are the first $m$ moments (or the first $2m$ moments for the CF transform).

According to the above relationships, the CF, MGF and PGF transform methods are essentially equivalent to the M method for values of $t_1,\ldots,t_k$ very close to zero (or very close to one for the PGF transform). In general, such choice of values is not optimal and a larger $k$ may be required for the M method to achieve the same level of efficiency as optimal, or nearly optimal, CF, MGF and PGF methods. This is because the moments may not carry as much information as the other transforms; a relevant heuristic argument is given in Kiefer [22]. We consider next a situation where the M method may be the most efficient.

PROPOSITION 2. *For* $k=1$ *and kernels* $g_t(Y)=\exp(itY)$ *and* $g_t(Y)=\exp(tY)$, *the necessary condition that*

$$(22) \qquad \max_t I_{G_n^*(t)}=\lim_{t\to 0}I_{G_n^*(t)} \qquad (=I_{G_n^*(1)})$$

*is*

$$(23) \quad \frac{\partial}{\partial\theta}E(Y_j^2|\mathcal{F}_{j-1})\operatorname{Var}(Y_j|\mathcal{F}_{j-1})-\frac{\partial}{\partial\theta}E(Y_j|\mathcal{F}_{j-1})\operatorname{Cov}(Y_j,Y_j^2|\mathcal{F}_{j-1})=0,$$

*for all $j$. An equivalent condition is*

$$(24) \qquad \frac{\partial}{\partial \theta} \operatorname{Var}(Y_j | \mathcal{F}_{j-1}) - \lambda_{j1} \frac{\partial}{\partial \theta} E(Y_j | \mathcal{F}_{j-1}) \operatorname{Var}^{1/2}(Y_j | \mathcal{F}_{j-1}) = 0,$$

*where $\lambda_{j1} = E[(Y_j - E(Y_j | \mathcal{F}_{j-1}))^3 | \mathcal{F}_{j-1}] / \operatorname{Var}^{3/2}(Y_j | \mathcal{F}_{j-1})$ is the index of skewness for the conditional distribution of $Y_j$. An analogous result holds for the kernel function $g_t(Y) = t^Y$ when $t$ approaches 1.*

PROOF. We obtain, with the help of a symbolic computation package,

$$\lim_{t \to 0} \frac{d}{dt} I_{G_n^*(t)} = \sum_{j=1}^{n} \frac{\frac{\partial}{\partial \theta} E(Y_j | \mathcal{F}_{j-1})}{\operatorname{Var}^2(Y_j | \mathcal{F}_{j-1})} \left[ \frac{\partial}{\partial \theta} E(Y_j^2 | \mathcal{F}_{j-1}) \operatorname{Var}(Y_j | \mathcal{F}_{j-1}) \right.$$
$$\left. - \frac{\partial}{\partial \theta} E(Y_j | \mathcal{F}_{j-1}) \operatorname{Cov}(Y_j, Y_j^2 | \mathcal{F}_{j-1}) \right].$$

Note that this limit is the coefficient of $t$ in the expansion of $I_{G_n^*(t)}$ in (21). Now set $\lim_{t \to 0} \frac{d}{dt} I_{G_n^*(t)} = 0$. This is equivalent to (23). Then the result follows on checking the sign of the second derivative. The equivalence of (23) and (24) is easy to prove. $\square$

Thus, the condition (23) is necessary for the quasi-score function $G_n^*(1)$, based on the first moment, to be more informative than the CF, MGF and PGF quasi-score functions based on a single point $t$. When $\operatorname{Var}(Y_j | \mathcal{F}_{j-1})$ is independent of $\theta$ for all $j$, as in stationary autoregressive processes, it follows from (24) that a necessary condition for (22) to hold is that $\lambda_{j1} = 0$ for all $j$. This occurs when the conditional distribution of $Y_j$ is symmetric around its mean. Although the form of this distribution is supposed to be unknown, its symmetry can be easily checked—the characteristic function of $Y_j - E(Y_j | \mathcal{F}_{j-1})$ must be real and even (see Rao [31], page 142).

5.2. *Choosing transform quasi-score functions.* Estimation using the estimating function $G_n^*(1)$ is the martingale version of what has been described by Wedderburn [36] as quasi-likelihood estimation in the context of independent observations. Notably, $G_n^*(1)$ is the exact score function for exponential family distributions with linear sufficient statistic (Remark 7, Section 3.2). Wedderburn's quasi-likelihood has been discussed in the context of discrete semimartingales in Hutton and Nelson [21], Godambe and Heyde [15] and Sørensen [35]. It should be noted that when $\sigma_j^2 = \operatorname{Var}(Y_j | \mathcal{F}_{j-1})$ is independent of $\theta$, $G_n^*(1)$ is the same as the weighted conditional least squares estimating function obtained by minimizing the sum of squares $\sum_{j=1}^{n} [Y_j - E(Y_j | \mathcal{F}_{j-1})]^2 / \sigma_j^2$ with respect to $\theta$. When $\sigma_j^2$ is also constant over all $j$, $G_n^*(1)$ reduces to the conditional least squares estimating function of Klimko and Nelson [23].

Although $G_n^*(1)$ will be globally optimal within the class of linear estimating functions based on the semimartingale representation of a stochastic process in an exponential family setting, it will typically be suboptimal within a larger class of estimating functions for nonexponential families. Then superior estimating functions may be readily constructed based on the semimartingale representation (5) of the transformed process.

First, consider polynomial estimating functions. A quadratic quasi-score function, essentially equivalent to $G_n^*(1, 2)$, has been considered by Godambe and Thompson [16] as an extension of the ordinary quasi-score $G_n^*(1)$ incorporating possible knowledge of the skewness and kurtosis of the underlying distribution. The formulation in Godambe and Thompson [16] involves a more general conditioning than the martingale structure, and the two components of the quadratic estimating function are the first two central moments corrected to zero mean and orthogonalized. A version of $G_n^*(1, 2)$ for i.i.d. variables had been considered earlier by Crowder [7].

Recall now from Remark 7 in Section 3.2 that (23) is the condition for the coefficient $w_{j2}$ in the quasi-score function $G_n^*(1, 2)$ to be zero. Therefore, $G_n^*(1)$ is as efficient as $G_n^*(1, 2)$ if the condition (23) holds for all $j$. Otherwise, when $\mathrm{Var}(Y_j | \mathcal{F}_{j-1})$ is independent of the parameter $\theta$, it is not difficult to show that $G_n^*(1, 2)$ reduces to $G_n^*(1)$ if the skewness is zero, and that the quasi-score function $G_n^*(1, 2, 3)$ also reduces to $G_n^*(1)$ if both skewness and kurtosis are zero. In situations different from those mentioned, $G_n^*(1)$ is less efficient than $G_n^*(1, 2)$. In fact, it is less efficient than a simple nonlinear transform quasi-score function [perturbed version of $G_n^*(1)$] that is based on a single optimal point $t$. Of course, only nonpolynomial martingale quasi-score functions are applicable if the variables $Y_j$ have no finite (conditional) moments, for example, if $Y_j$ has a stable distribution; see Example 1 in the next section.

A composite transform quasi-score function utilizing more than a single point $t$ will generally be more efficient. We may also combine different transforms. For example, we may choose $\mathbf{g}(Y_j) = (Y_j, Y_j^2, \ldots, Y_j^{k_1}, \exp(t_1 Y_j), \ldots, \exp(t_{k_2} Y_j))'$, with optimal points $t_{k_1+1}, \ldots, t_{k_2}$, thereby combining the convenience of the M method with the higher efficiency of the MGF method.

Composite quasi-score functions may be also effective in dealing with situations in which the estimating functions for a vector parameter in alternative methods (e.g., conditional least squares) are functionally dependent and, thus, no estimates can be obtained; see Example 2 in the next section.

**6. Examples.** The following two brief examples illustrate important features of the proposed method. A more detailed study of these applications will be reported elsewhere.

EXAMPLE 1. *AR(1) process with symmetric stable error.* Consider the AR(1) process

$$Y_j = \phi Y_{j-1} + \varepsilon_j,$$

where $\{\varepsilon_j\}$ is an i.i.d. sequence of symmetric stable random variables with characteristic function $c(t) = \exp(-|t|^\alpha)$, $0 < \alpha \le 2$. Closed form density representation exists only in the Cauchy $(\alpha = 1)$ and Gaussian $(\alpha = 2)$ cases. We wish to estimate $\phi$ on the basis of a sample $Y_1, \ldots, Y_n$. Consider the martingale difference

$$h_j(t) = \exp(itY_j) - E(\exp(itY_j)|Y_{j-1}) = \exp(itY_j) - \exp(it\phi Y_{j-1} - |t|^\alpha).$$

Then the transform quasi-score function constructed using the real kernel $(\cos(tY_j), \sin(tY_j))$ is $G_n^*(t) = \sum_{j=1}^n \mathbf{w}_j^* \mathbf{h}_j$, where $\mathbf{h}_j = (\operatorname{Re} h_j(t), \operatorname{Im} h_j(t))'$ and $\mathbf{w}_j^* = (E(\dot{\mathbf{h}}_j|Y_{j-1}))'(E(\mathbf{h}_j \mathbf{h}_j'|Y_{j-1}))^{-1}$. Noticing that $(E(\dot{\mathbf{h}}_j|Y_{j-1}))' = (\frac{\partial}{\partial \phi}\operatorname{Re} h_j(t), \frac{\partial}{\partial \phi}\operatorname{Im} h_j(t))'$, and using properties of the cosine and sine functions to derive the entries of $(E(\mathbf{h}_j \mathbf{h}_j'|Y_{j-1}))^{-1}$, we can show that

$$G_n^*(t) = \frac{2t\exp(2^\alpha|t|^\alpha)}{\exp(|t|^\alpha)(\exp(2^\alpha|t|^\alpha) - 1)} \sum_{j=1}^n Y_{j-1}\sin(t(\phi Y_{j-1} - Y_j))$$

and

$$(25) \qquad I_{G_n^*(t)} = \frac{2t^2\exp(2^\alpha|t|^\alpha)}{\exp(2|t|^\alpha)(\exp(2^\alpha|t|^\alpha) - 1)} \sum_{j=1}^n Y_{j-1}^2.$$

It is important to note that because of the infinite variance of the $Y_j$'s, the $O_F$ information quantity $E(G_n^{*2}(t))$, which is equal to $E(I_{G_n^*(t)})$ in the finite variance case, is not defined in the present case, implying that the $O_F$-optimality criterion is not applicable. However, $I_{G_n^*(t)} < \infty$ for all $n \ge 1$ and given $t \in T$.

Turning now to the choice of the optimal value of the point $t \in T$, we observe that this choice is independent of the parameter $\phi$ and the observations. We also observe in (25) that the factor multiplying $\sum_{j=1}^n Y_{j-1}^2$ in $I_{G_n^*(t)}$ is an even function of $t$. We consider this factor for $t \ge 0$ to obtain its generalized series expansion $2^{1-\alpha}t^{2-\alpha} + 2^{1-\alpha}(2^{\alpha-1} - 2)t^2 + o(t^{\alpha+2})$. It follows that

$$\lim_{t \to 0} I_{G_n^*(t)} = \begin{cases} 0, & \text{if } 0 < \alpha < 2, \\ 1/2 \sum_{j=1}^n Y_{j-1}^2, & \text{if } \alpha = 2. \end{cases}$$

In the case of a normal distribution $(\alpha = 2)$, $I_{G_n^*(t)} < 1/2\sum_{j=1}^n Y_{j-1}^2 = \lim_{t \to 0} I_{G_n^*(t)}$ $(= I_{S_n})$, that is, $I_{G_n^*(t)}$ attains its maximum at the origin [in accordance with (22)]. For $\alpha < 2$, $I_{G_n^*(t)}$ attains its maximum at a value of $t > 0$; the smaller the value of $\alpha$, the larger the optimal point $t^*$.

TABLE 1
*The efficiency of $G_n^*(t^*)$*

| $\alpha$ | Fisher information (i.i.d.) | $I_{G_n^*(t^*)}/\sum_{j=1}^{n}Y_{j-1}^2$ | $t^*$ | eff$_c\left(G_n^*(t^*),S_n\right)$ |
|---|---|---|---|---|
| 2.0 | 0.500 | 0.500 | 0.00 | 1 |
| 1.9 | 0.473 | 0.469 | 0.3852 | 0.991 |
| 1.7 | 0.442 | 0.428 | 0.4767 | 0.968 |
| 1.5 | 0.428 | 0.391 | 0.5384 | 0.913 |
| 1.3 | 0.431 | 0.358 | 0.6087 | 0.831 |
| 1.1 | 0.463 | 0.332 | 0.7148 | 0.717 |
| 1.0 | 0.500 | 0.324 | 0.7968 | 0.648 |
| 0.8 | 0.678 | 0.330 | 1.1022 | 0.487 |

We can assess the efficiency of $G_n^*(t^*)$ for varying $\alpha$ by comparing the factor multiplying $\sum_{j=1}^{n}Y_{j-1}^2$ in $I_{G_n^*(t^*)}$ with the essentially exact Fisher information (per observation) for the location parameter, say, $\phi$, computed by DuMouchel [8] in the setting of i.i.d. variables for selected values of $\alpha$. This is presented in Table 1. The optimal values for $t$ in each case are also shown.

It can be seen that for the range of values of $\alpha$ reported in Table 1 the efficiency of $G_n^*(t^*)$ decreases as $\alpha$ deviates from $\alpha = 2$. The efficiency of the transform quasi-score function can be increased by introducing more points $t$, albeit at the expense of increased computational complexity.

EXAMPLE 2. *A first-order gamma autoregressive model.* Consider the first-order gamma autoregressive model (Sim [33]) given by

$$(26) \qquad Y_j = \alpha * Y_{j-1} + \varepsilon_j,$$

where the operator $*$ is defined as $\alpha * Y = \sum_{i=1}^{N(Y)} E_i$, and (i) the $\varepsilon_j$ are i.i.d. Gamma$(\alpha, \nu)$ random variables with $\alpha,\ \nu > 0$; (ii) the $E_i$ are i.i.d. exponential $(\alpha)$ random variables; (iii) for each fixed positive value of $y$, $N(y)$ is a Poisson random variable with parameter $\lambda = p\alpha$, and $0 \le p < 1$.

Expression (26) is an autoregressive representation for a stationary gamma process whose joint density function has a certain type of Laplace transform. This process has been used in the study of stochastic reservoir systems with Markovian inflows; see Sim [33] and references therein. The conditional Laplace transform of $Y_j$ can be easily derived (Sim [33]) and is given by

$$(27) \qquad E(\exp(-sY_j)|Y_{j-1} = y_{j-1}) = \left(\frac{\alpha}{\alpha + s}\right)^{\nu}\exp\left(-\frac{\lambda s y_{j-1}}{\alpha + s}\right).$$

Inversion of (27) gives the conditional density of $Y_j$ as

$$f(y_j|y_{j-1}) = \alpha \left( \frac{\alpha y_j}{\lambda y_{j-1}} \right)^{(\nu-1)/2} \exp[-(\alpha y_j + \lambda y_{j-1})] I_{\nu-1}[2(\lambda \alpha y_j y_{j-1})^{1/2}],$$

where $I_r(z)$ is the modified Bessel function of the first kind and of order $r$. The likelihood function being complicated, Sim suggested (but did not apply) the conditional least squares method for the estimation of the three parameters of the model. The conditional least squares equations are

$$\sum_{j=1}^{n} Y_{j-1}[\alpha Y_j - \lambda Y_{j-1} - \nu] = 0,$$

$$\sum_{j=1}^{n} (\lambda Y_{j-1} + \nu)[\alpha Y_j - \lambda Y_{j-1} - \nu] = 0,$$

$$\sum_{j=1}^{n} [\alpha Y_j - \lambda Y_{j-1} - \nu] = 0.$$

This system of equations has no solution other than the trivial solution $\lambda = \alpha = \nu = 0$. The same is true for the quasi-likelihood equations $\mathbf{G}_n^*(1) = \mathbf{0}$, and for the quasi-likelihood equations $\mathbf{G}_n^*(s) = \mathbf{0}$ based on the Laplace transform (27) and employing a single point $s$. All three systems of equations mentioned above have a solution if one of the parameters is fixed. In particular, if we consider the parameter $\nu$ fixed, then maximum likelihood estimation is also possible, though very complicated.

Estimation of all three parameters of the model is possible by using a composite quasi-score function based on the Laplace transform (27). Conditional moments of any order may be obtained from (27) in simple form; see Sim [33]. Then the convenient composite quasi-score function $\mathbf{G}_n^*(1,2)$ can be used. The corresponding system of quasi-likelihood equations, shown in simulations to have a numerical solution, is

$$\sum_{j=1}^{n} [5\lambda^3 Y_{j-1}^4 + (-6\alpha Y_j + 10\nu + 4)\lambda^2 Y_{j-1}^3$$

$$+ (\alpha^2 Y_j^2 - 6((1+\nu)\alpha Y_j - \nu - \nu^2))\lambda Y_{j-1}^2$$

$$+ ((-\nu - \nu^2)\alpha Y_j + \nu^2 + \nu^3)Y_{j-1}] = 0,$$

$$\sum_{j=1}^{n} [-3\lambda^4 Y_{j-1}^4 + (2\alpha Y_j - 10\nu - 2)\lambda^3 Y_{j-1}^3$$

$$+ (\alpha^2 Y_j^2 + 6\nu\alpha Y_j - 6\nu - 12\nu^2)\lambda^2 Y_{j-1}^2$$

$$- (6\nu - 5\alpha Y_j)\nu(1+\nu)\lambda Y_{j-1} + (\nu^2 + \nu^3)\alpha Y_j - \nu^3 - \nu^4] = 0,$$

$$\sum_{j=1}^{n}[6\lambda^3 Y_{j-1}^3 + ((25/2)\nu + 5 - 8\alpha Y_j)\lambda^2 Y_{j-1}^2$$

$$+ (2\alpha^2 Y_j^2 - 9(1+\nu)\alpha Y_j + 8\nu + 8\nu^2)\lambda Y_{j-1}$$

$$+ (1/2)\nu Y_j^2 \alpha^2 - 2\nu(1+\nu)\alpha Y_j + (3/2)\nu^2(1+\nu)] = 0.$$

**Acknowledgments.** I wish to thank Professor Mary E. Thompson for introducing me to this research area and for inspiring me in an earlier phase of my work on the subject. I also wish to thank the referees for their comments which have led to a substantial improvement of the paper.

## REFERENCES

[1] ABRAHAM, B. and BALAKRISHNA, N. (1999). Inverse Gaussian autoregressive models. *J. Time Ser. Anal.* **20** 605–618. MR1749577

[2] ALZAID, A. A. and AL-OSH, M. A. (1990). An integer-valued *p*th-order autoregressive structure (INAR(*p*)) process. *J. Appl. Probab.* **27** 314–324. MR1052303

[3] BILLARD, L. and MOHAMED, F. Y. (1991). Estimation of the parameters of an EAR(*p*) process. *J. Time Ser. Anal.* **12** 179–192. MR1128943

[4] BRANT, R. (1984). Approximate likelihood and probability calculations based on transforms. *Ann. Statist.* **12** 989–1005. MR0751287

[5] BROCKWELL, P. J. and LIU, J. (1992). Estimating the noise parameters from observations of a linear process with stable innovations. *J. Statist. Plann. Inference* **33** 175–186. MR1190618

[6] BURRILL, C. W. (1972). *Measure, Integration and Probability.* McGraw-Hill, New York. MR0457657

[7] CROWDER, M. (1987). On linear and quadratic estimating functions. *Biometrika* **74** 591–597. MR0909363

[8] DUMOUCHEL, W. H. (1975). Stable distributions in statistical inference. II. Information from stably distributed samples. *J. Amer. Statist. Assoc.* **70** 386–393. MR0378191

[9] EPPS, T. W. and PULLEY, L. B. (1985). Parameter estimates and tests of fit for infinite mixture distributions. *Comm. Statist. Theory Methods* **14** 3125–3145.

[10] FEIGIN, P. D., TWEEDIE, R. L. and BELYEA, C. (1983). Weighted area techniques for explicit parameter estimation in multistage models. *Austral. J. Statist.* **25** 1–16. MR0712001

[11] FEUERVERGER, A. (1990). An efficiency result for the empirical characteristic function in stationary time-series models. *Canad. J. Statist.* **18** 155–161. MR1067167

[12] FEUERVERGER, A. and McDUNNOUGH, P. (1981). On some Fourier methods for inference. *J. Amer. Statist. Assoc.* **76** 379–387. MR0624339

[13] FEUERVERGER, A. and McDUNNOUGH, P. (1981). On the efficiency of empirical characteristic function procedures. *J. Roy. Statist. Soc. Ser. B* **43** 20–27. MR0610372

[14] FEUERVERGER, A. and McDUNNOUGH, P. (1984). On statistical transform methods and their efficiency. *Canad. J. Statist.* **12** 303–317. MR0782952

[15] GODAMBE, V. P. and HEYDE, C. C. (1987). Quasi-likelihood and optimal estimation. *Internat. Statist. Rev.* **55** 231–274. MR0963141

[16] GODAMBE, V. P. and THOMPSON, M. E. (1989). An extension of quasi-likelihood estimation (with discussion). *J. Statist. Plann. Inference* **22** 137–172. MR1004344

[17] GRUNWALD, G. K., HYNDMAN, R. J., TEDESCO, L. and TWEEDIE, R. L. (2000). Non-Gaussian conditional linear AR(1) models. *Aust. N. Z. J. Stat.* **42** 479–495. MR1802969

[18] HEYDE, C. C. (1987). On combining quasilikelihood estimating functions. *Stochastic Process. Appl.* **25** 281–287. MR0915142

[19] HEYDE, C. C. (1997). *Quasi-Likelihood and Its Application: A General Approach to Optimal Parameter Estimation.* Springer, New York. MR1461808

[20] HOETING, J. A., TWEEDIE, R. L. and OLVER, C. S. (2003). Transform estimation of parameters for stage-frequency data. *J. Amer. Statist. Assoc.* **98** 503–514. MR2011667

[21] HUTTON, J. E. and NELSON, P. I. (1986). Quasilikelihood estimation for semimartingales. *Stochastic Process. Appl.* **22** 245–257. MR0860935

[22] KIEFER, N. M. (1978). Comment on "Estimating mixtures of normal distributions and switching regressions," by R. E. Quandt and J. B. Ramsay. *J. Amer. Statist. Assoc.* **73** 744–745.

[23] KLIMKO, L. A. and NELSON, P. I. (1978). On conditional least squares estimation for stochastic processes. *Ann. Statist.* **6** 629–642. MR0494770

[24] LEITNAKER, M. G. (1989). Estimation of delay times in stochastic compartmental models. *Biometrics* **45** 1239–1247.

[25] MCCULLOCH, J. H. (1996). Financial applications of stable distributions. In *Statistical Methods in Finance* (G. S. Maddala and C. R. Rao, eds.) 393–425. North-Holland, Amsterdam. MR1602156

[26] MCLEISH, D. L. (1984). Estimation for aggregate models: The aggregate Markov chain (with discussion). *Canad. J. Statist.* **12** 265–285. MR0782050

[27] MERKOURIS, T. (1992). A transform method for optimal estimation in stochastic processes: Basic aspects. In *Proc. Symposium in Honour of Professor V. P. Godambe* (J. Chen, ed.). Univ. Waterloo, Waterloo, Canada.

[28] MERKOURIS, T. (1992). A transform method for optimal estimation in stochastic processes. Ph.D. dissertation, Dept. Statistics, Univ. Waterloo.

[29] MURPHY, S. and LI, B. (1995). Projected partial likelihood and its application to longitudinal data. *Biometrika* **82** 399–406. MR1354237

[30] NIKIAS, C. L. and SHAO, M. (1995). *Signal Processing with Alpha-Stable Distributions and Applications.* Wiley, New York.

[31] RAO, C. R. (1973). *Linear Statistical Inference and Its Applications*, 2nd ed. Wiley, New York. MR0346957

[32] SCHUH, H.-J. and TWEEDIE, R. L. (1979). Parameter estimation using transform estimation in time-evolving models. *Math. Biosci.* **45** 37–67. MR0535954

[33] SIM, C. H. (1990). First-order autoregressive models for gamma and exponential processes. *J. Appl. Probab.* **27** 325–332. MR1052304

[34] SMALL, C. G. and MCLEISH, D. L. (1994). *Hilbert Space Methods in Probability and Statistical Inference.* Wiley, New York. MR1269321

[35] SØRENSEN, M. (1990). On quasi-likelihood for semimartingales. *Stochastic Process. Appl.* **35** 331–346. MR1067117

[36] WEDDERBURN, R. W. M. (1974). Quasi-likelihood functions, generalized linear models and the Gauss–Newton method. *Biometrika* **61** 439–447. MR0375592

[37] YAO, Q. and MORGAN, B. J. T. (1999). Empirical transform estimation for indexed stochastic models. *J. R. Stat. Soc. Ser. B Stat. Methodol.* **61** 127–141. MR1664112

Statistics Canada
R. H. Coats Bldg., 16-D
Tunney's Pasture
Ottawa, Ontario
Canada K1A 0T6
E-mail: Takis.Merkouris@statcan.ca